# COEXISTENCE FOR A MULTITYPE CONTACT PROCESS WITH SEASONS

By B. Chan, R. Durrett and N. Lanchier

*Cornell University, Cornell University and Arizona State University*

We introduce a multitype contact process with temporal heterogeneity involving two species competing for space on the $d$-dimensional integer lattice. Time is divided into seasons called alternately season 1 and season 2. We prove that there is an open set of the parameters for which both species can coexist when their dispersal range is large enough. Numerical simulations also suggest that three species can coexist in the presence of two seasons. This contrasts with the long-term behavior of the time-homogeneous multitype contact process for which the species with the higher birth rate outcompetes the other species when the death rates are equal.

**1. Introduction.** Consider an ecosystem inhabited by multiple species competing for a single resource represented by a function $R$. Let $u_i$, $i = 1, 2, \ldots, K$, denote the population density of species $i$. In the original models of Lotka (1925) and Volterra (1928), the per capita growth rate for each species increases linearly with the amount of resource available:

$$(1) \qquad \frac{1}{u_i} \frac{du_i}{dt} = \beta_i R(u_1, u_2, \ldots, u_K) - \delta_i,$$

where $\beta_i$ and $\delta_i$, respectively, denote the birth and death rates of species $i$. When $R$ is a decreasing function of population densities, Volterra (1928) showed that only the species with the highest $\beta_i$ to $\delta_i$ ratio survives as $t \to \infty$.

When $R = 1 - (u_1 + u_2 + \cdots + u_K)$ denotes the density of space available, the Lotka–Volterra model (1) can be seen as the mean-field approximation, or nonspatial deterministic analogue, of the multitype contact process introduced by Neuhauser (1992). Assuming for simplicity that the number of species is $K = 2$, her model is a continuous-time Markov process whose state at time $t$ is a function $\eta_t$ that maps the $d$-dimensional integer lattice $\mathbb{Z}^d$ into $\{0, 1, 2\}$, with 0 denoting empty sites, and 1 and 2 denoting sites occupied by an individual of species 1 and 2, respectively. The evolution at site $x \in \mathbb{Z}^d$ is described by the following transition rates:

$$0 \to 1 \qquad \text{at rate } \beta_1 \sum_{y \in \mathbb{Z}^d} p(y, x) \mathbb{1}_{\{\eta(y) = 1\}}, \qquad 1 \to 0 \qquad \text{at rate } \delta_1,$$

$$0 \to 2 \qquad \text{at rate } \beta_2 \sum_{y \in \mathbb{Z}^d} p(y, x) \mathbb{1}_{\{\eta(y) = 2\}}, \qquad 2 \to 0 \qquad \text{at rate } \delta_2,$$

where $p(y, x)$ is a translation invariant transition probability on $\mathbb{Z}^d$. This indicates that species of type $i$, $i = 1, 2$, produces offspring at rate $\beta_i$ and dies at rate $\delta_i$. Moreover, an offspring produced at site $y \in \mathbb{Z}^d$ is sent to a site $x$ chosen at random with probability $p(y, x)$. If site $x$ is vacant, it becomes occupied by the offspring, otherwise the birth is suppressed. Neuhauser (1992) proved for the multitype contact process that if $\delta_1 = \delta_2$ then the species with the larger birth rate outcompetes the other species. She also conjectured that, as predicted by the model (1), the species with the larger $\beta_i$ to $\delta_i$ ratio takes over. While it is believed that there is no coexistence in the Neuhauser's multitype contact process, there are a number of interacting particle systems for which coexistence has been proved rigorously.

One mechanism that allows coexistence is when intraspecific competition exceeds interspecific competition. This has been proved by Neuhauser and Pacala (1999) for a voter type model in which each site of the lattice is occupied by an individual of type 1 or 2. When both species have the same birth rate, the state at site $x \in \mathbb{Z}^d$ evolves according to the transition rates

$$1 \to 2 \qquad \text{at rate } f_2(f_1 + \alpha_{12} f_2),$$

$$2 \to 1 \qquad \text{at rate } f_1(f_2 + \alpha_{21} f_1),$$

where $f_i = f_i(x)$ denotes the fraction of neighbors of site $x$ in state $i$. The parameters $\alpha_{ij}$ model the pressure of species $i$ on species $j$. Except for the 1-dimensional nearest neighbor case, in the absence of interspecific competition, that is, $\alpha_{12} = \alpha_{21} = 0$, there is a stationary distribution in which sites are independent [see Theorem 1 in Neuhauser and Pacala (1999)], indicating that both species coexist. While increasing $\alpha_{12}$ and $\alpha_{21}$, spatial correlations build up, but coexistence is maintained as long as the interspecific competition is sufficiently small. The inclusion of a spatial structure, however, translates into a reduction of the set of parameters for which coexistence occurs.

A second mechanism that promotes coexistence is environmental heterogeneity. Suppose that the regular lattice $\mathbb{Z}^2$ is inhabited by two types of plants. The environment is static and described by a partition $\{H_1, H_2\}$ of the space with $H_i$ denoting the set of sites occupied by plant $i$. Consider two specialist symbionts; for $i = 1, 2$, symbiont of type $i$ feeds on plant $i$. This can be modeled by a multitype contact process in which births of symbionts of type $i$ are suppressed whenever the offspring is sent to a site $x \in H_j$ with $i \neq j$. Since the habitats $H_1$ and $H_2$ are disjoint and so symbionts 1 and 2 do not interact, the problem of coexistence reduces to the following percolation problem: If for $i = 1, 2$, there exists an infinite sequence of sites in $H_i$ along which symbiont $i$ can spread then coexistence is possible. This depends on both the spatial arrangement of $H_i$ and the dispersal range of symbiont $i$.

Lanchier and Neuhauser (2006) proved that coexistence occurs in a non-trivial case in which the habitats of both symbionts overlap. Their model is referred to as the static-host model. Symbiont 1 is a specialist feeding on plant 1 while symbiont 2 is a generalist feeding on both plants. The universe is tiled with $N \times N$ squares. Whether a square is occupied by plants of type 1 or 2 is determined by flipping independent fair coins. In this situation, coexistence occurs when the birth rate of the specialist is larger than that of the generalist and the square size $N$ is much larger than the dispersal range of the symbionts. In the mean-field approximation of the static-host model, the birth rate of the specialist needs to exceed twice the birth rate of the generalist so that coexistence occurs. Lanchier and Neuhauser (2006) also introduced a stochastic process in which spatial heterogeneities are no longer static but generated by the symbionts. This model is referred to as the dynamic-host model. In the absence of symbionts, plants of type 1 and 2 compete according to a biased voter model with a selective advantage for the 1's. In particular, plants of type 1 outcompete plants of type 2. Introducing a specialist pathogen feeding on plants of type 1 whose harmful effect is modeled by an increase of the death rate of infected plants, Durrett and Lanchier (2008) proved that coexistence may occur provided the dispersal range of the species is sufficiently large, but conjectured that coexistence is not possible if the symbiont is a specialist mutualist feeding on plants of type 2.

In some cases, the presence of two types of space is not immediately obvious. Chan and Durrett (2006) considered a generalization of the multitype contact process in which species 1 is a good competitor with a short range dispersal kernel, whereas species 2 has a long range dispersal kernel which makes it a good colonizer. In this situation, occurrences of forest fires (i.e., removal of all the individuals contained in a large square) at an appropriate rate allows species 2 to survive by migrating to newly created patches. In this model, there are two types of space: blocks freshly cleared by a fire which are occupied by individuals of species 2 and older blocks which are occupied by individuals of species 1.

The aim of this article is to prove that temporal heterogeneity can also promote coexistence in a stochastic spatial model. Armstrong and McGehee (1976) have shown in an ODE model that coexistence can occur for species living in a periodic environment with seasonal changes in birth rates: when each species has a specific growing season disjoint from all others, and all resources are conservative (meaning it is a linear function of population densities), a system with four resources can support more than four species without any one going into extinction. Our objective is to extend their observation to competition models based on interacting particle systems.

Recall that when two contact processes with equal death rates $\delta$ and equal dispersal kernels compete with one another on $\mathbb{Z}^2$, one of the two species outcompetes the other one. To search for coexistence, we follow the work of Armstrong and McGehee and consider an environment with seasonal changes in birth rates: time will be divided into seasons of length $D$, called alternately season 1 and season 2. For simplicity and because this is the case of interest in spatial ecology, we assume that $d = 2$. However, all our results extend easily to any dimension $d \geq 1$. Since we will deal with long range interactions, it is convenient to assume that particles evolve on the rescaled lattice $\mathbb{Z}^2/L$ where $L$ is a large integer which is referred to as the range of the interactions. Our model is a continuous-time Markov process whose state at time $t$ is a function $\xi_t$ that maps the rescaled lattice into $\{0, 1, 2\}$. To formulate the dynamics, we let

$$\mathcal{N}(x) = \{y \in \mathbb{Z}^2/L : 0 < \|x - y\|_\infty \leq 1\}$$

be the interaction neighborhood of site $x \in \mathbb{Z}^2/L$, and

$$f_i(x, \xi) = \frac{|\{y \in \mathcal{N}(x) : \xi(y) = i\}|}{|\mathcal{N}(x)|}$$

be the fraction of neighbors of site $x$ in state $i$. The state of site $x \in \mathbb{Z}^2/L$ flips according to the following transition rates:

| transition | season 1 | season 2 |
|:---:|:---:|:---:|
| $0 \to 1$ | $\beta_{11} f_1(x, \xi)$ | $\beta_{12} f_1(x, \xi)$ |
| $0 \to 2$ | $\beta_{21} f_2(x, \xi)$ | $\beta_{22} f_2(x, \xi)$ |
| $1 \to 0$ | $\delta_1$ | $\delta_1$ |
| $2 \to 0$ | $\delta_2$ | $\delta_2$ |

This indicates that each species behaves individually like a contact process. Their associated birth rates, however, depend on the season, making our process a generalization of the Neuhauser's competing model in which the birth rates are functions of time. We call the process $\xi_t$ the *periodic competition model*. We assume for simplicity that each of the species has the same death rate in both seasons. However, the reduction to equal death rates is not important for the proofs. In particular, Theorem 1 and Corollary 3 below can be extended to the case when the death rates may vary provided one replaces the death rate by the average of the death rates in season 1 and season 2, respectively. We now give a sufficient condition for coexistence.

To begin with, we observe that when the dispersal range $L$ tends to infinity, the population dynamics (in a finite volume) converges to the following mean-field model:

$$
(2) \qquad \begin{aligned}
\frac{du_1}{dt} &= \beta^1(t)u_1(1 - u_1 - u_2) - \delta_1 u_1, \\
\frac{du_2}{dt} &= \beta^2(t)u_2(1 - u_1 - u_2) - \delta_2 u_2,
\end{aligned}
$$

where $u_i$, $i = 1, 2$, denotes the population density of species $i$ and $\beta^i(t)$ are periodic step functions with period $2D$ given by

$$
\beta^i(t) = \begin{cases} \beta_{i1}, & \text{when } 0 \le t < D, \\ \beta_{i2}, & \text{when } D \le t < 2D. \end{cases}
$$

In the absence of species 1, species 2 evolves according to the pair of logistic equations

$$
(3) \qquad \frac{du_2}{dt} = \begin{cases} \beta_{21}u_2(1 - u_2) - \delta_2 u_2, & \text{if } t \in [2nD, (2n+1)D), \\ \beta_{22}u_2(1 - u_2) - \delta_2 u_2, & \text{if } t \in [(2n+1)D, (2n+2)D). \end{cases}
$$

Denoting by $u_2(t)$ the solution, the function $t \mapsto u_2(2nD + t)$ for $0 \le t \le 2D$ converges uniformly to a piecewise-smooth periodic curve $\bar{u}_2$. This follows from the monotonicity of the solution with respect to the initial condition. Moreover, since $\bar{u}_2$ is periodic and continuous, standard topological arguments imply that the convergence is uniform on the real line. We call $\bar{u}_2$ the *equilibrium curve for the 2's.* One can define $\bar{u}_1$ for the 1's in a similar way. The top panel of Figure 1 gives an illustration of the density of 2's in the absence of 1's for the interacting particle system which, in view of the size of the universe and the range of the interactions, is a good approximation of $\bar{u}_2$.

We now return to the spatial model. When the birth rates $\beta_{ij}$ are chosen appropriately, one can compare the 1's (resp., the 2's) with a supercritical branching process and show that the 1's and 2's coexist. Based on Armstrong and McGehee ([1976](#)), we expect that any number of species can coexist provided there are as many seasons as there are species, but for simplicity we restrict ourselves to a system with only two species.

THEOREM 1. *Suppose that*

$$
\frac{1}{2D} \int_0^{2D} \beta^1(t)(1 - \bar{u}_2(t))\, dt > \delta_1 \quad \text{and} \quad \frac{1}{2D} \int_0^{2D} \beta^2(t)(1 - \bar{u}_1(t))\, dt > \delta_2.
$$

*Then, when $L$ is large, coexistence occurs.*

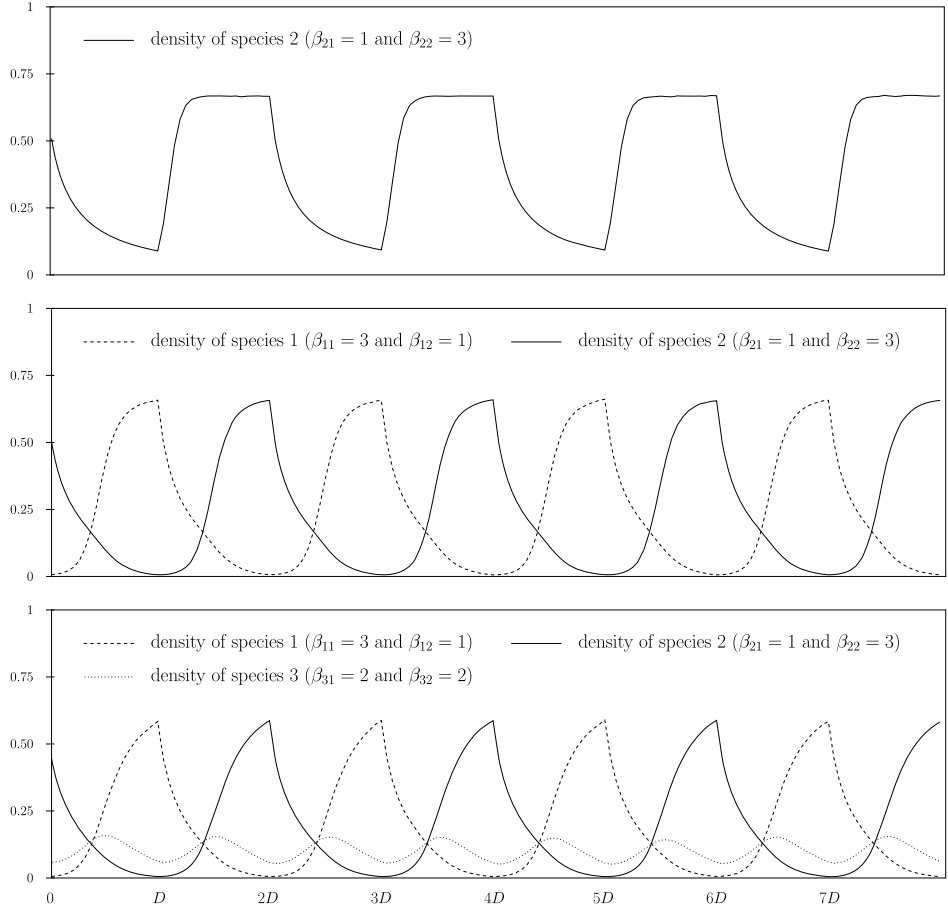

Fig. 1. *Evolution of the densities of species for the interacting particle system along 4 periods = 8 seasons. The universe is the $400 \times 400$ lattice with periodic boundary conditions. The length of each season is $D = 10$, the range of the interactions $L = 200$ and the death rates $\delta_i = 1$. Our pictures show species 2 alone, species 1 and 2, and species 1, 2 and 3 in competition, respectively.*

In the previous theorem and after, coexistence for the interacting particle system means that starting from any initial configuration with infinitely many 1's and infinitely many 2's,

$$\liminf_{t \to \infty} P(\xi_t(x) = i) > 0 \qquad \text{for } i = 1, 2 \text{ and all } x \in \mathbb{Z}^2 / L.$$

Note that due to temporal heterogeneities there is no stationary distribution, but there is an initial distribution so that the states at times 0 and $2D$ have the same distribution. It seems reasonable to call this a periodic distribution. Since $\int (1 - \bar{u}_2)/2D \, dt$ is the average number of empty sites available for invasion, the first integral in Theorem 1 represents the growth rate of 1's when rate averaged over time. Therefore, if the average growth rate is greater than the death rate then the 1's can invade a community of 2's. Similarly, the second inequality says that if the average growth rate for species 2 is greater than the death rate then the 2's can invade a community of 1's.

As previously mentioned, it has already been proved by Armstrong and McGehee (1976) that temporal heterogeneity promotes coexistence. Their result, however, was based on a deterministic nonspatial model. It is important to extend their result to a model including both temporal heterogeneity and spatial structure in the form of local interactions because the spatial component is identified in ecology as an important factor in how communities are shaped. Even though the sufficient conditions we found for coexistence of two species in the presence of two seasons turn out to be the same for both the interacting particle system and its nonspatial analog, our main result is interesting for ecologists since it is known from past research that spatial models can result in predictions that differ from nonspatial models. For instance, it has been proved theoretically in Neuhauser and Pacala (1999) and Lanchier and Neuhauser (2009) that the inclusion of local interactions may translate into a reduction of the coexistence region. Their models were based on particle systems with short range interactions. Interestingly, Neuhauser (1994) also proved that a sexual reproduction process with a long range dispersal kernel does not exhibit the same features as the corresponding mean-field model, namely, for intermediate values of the birth rate survival occurs for the mean-field model while the interacting particle system exhibits a metastable behavior and eventually goes extinct. In other words, even in the presence of long range interactions, spatial and nonspatial models may have significantly different behaviors.

The second panel of Figure 1 shows the evolution of the densities of species 1 and 2 in competition, which illustrates Theorem 1. The bottom panel suggests that temporal heterogeneity can also promote coexistence of three species even in an environment with only two seasons. Interestingly, numerical simulations also suggest that, with the birth rates indicated in Figure 1, the three species coexist whereas in the absence of species 2, species 3 outcompetes species 1. This indicates that the inclusion of an additional species may promote coexistence. In conclusion, we state the following.

CONJECTURE 2.  *In the presence of two seasons, three species can coexist.*

The conjecture is probably very difficult to prove. In principle it can be approached with invadability ideas as in Durrett (2002). Suppose for simplicity that we have found three species that can coexist in pairs. To show that all three can coexist, we need to (i) show that the two species systems converge to equilibrium, and that (ii) in all three cases the third species can invade the other two in their equilibrium. Problem (i) seems difficult. Indeed, we do not see how to prove that the mean field ODE converges to equilibrium. Due to the time inhomogeneity, one does not have Lyapunov functions.

We now derive from the two integrals in Theorem 1 a (stronger but) more explicit sufficient condition for the coexistence of both species in the spatially explicit model. This condition is given in Corollary 3 below. First of all, we set

$$p_{2j} = \lim_{n \to \infty} u_2(2nD + (j-1)D) \qquad \text{while holding } u_1 \equiv 0$$

for $j = 1, 2$. That is, $p_{2j}$ is the limiting density of the 2's at the beginning of season $j$ in the absence of 1's. Since the 2's grow according to the pair of logistic equations (3), the densities $p_{21}$ and $p_{22}$ are related by

$$\rho(p_{21}, D; \beta_{21}, \delta_2) = p_{22} \quad \text{and} \quad \rho(p_{22}, D; \beta_{22}, \delta_2) = p_{21},$$

where the function $\rho$ is defined by

$$\rho(u_0, t; \beta, \delta) = \frac{u_0 \exp(rt)}{1 - K^{-1} u_0 (1 - \exp(rt))},$$

where

$$r = \beta - \delta \quad \text{and} \quad K = 1 - \delta/\beta$$

are the intrinsic rate of growth and the carrying capacity, respectively. By letting $j = (n \bmod 2) + 1$, one can define $\bar{u}_2$ recursively, by

$$\bar{u}_2(0) = p_{21} \quad \text{and} \quad \bar{u}_2(nD + \theta) = \rho(\bar{u}_2(nD), \theta; \beta_{2j}, \delta_2)$$

for $0 \leq \theta < D$. One can define $p_{11}$, $p_{12}$ and $\bar{u}_1$ for the 1's in a similar way. Now, assuming that season 2 is the growing season for species 2, that is, $\beta_{21} < \beta_{22}$, the limiting curve $\bar{u}_2$ is convex and decreasing in season 1 and concave and increasing in season 2. (See the top panel of Figure 1 for an approximation of $\bar{u}_2$.) In particular,

$$\frac{1}{D} \int_0^D 1 - \bar{u}_2(t) \, dt \geq 1 - \frac{p_{21} + p_{22}}{2} \quad \text{and} \quad \frac{1}{D} \int_D^{2D} 1 - \bar{u}_2(t) \, dt \geq 1 - p_{21}.$$

This indicates that the densities of empty sites available for the 1's in seasons 1 and 2 are bounded from below by $1 - (p_{21} + p_{22})/2$ and $1 - p_{21}$, respectively. Therefore, the 1's dominate a contact process with birth rate

$$\beta_{11}\left(1 - \frac{p_{21} + p_{22}}{2}\right) \quad \text{and} \quad \beta_{12}(1 - p_{21})$$

in seasons 1 and 2, respectively. Similarly, assuming that $\beta_{11} > \beta_{12}$, the 2's dominate a contact process with birth rate

$$\beta_{21}(1 - p_{12}) \quad \text{and} \quad \beta_{22}\left(1 - \frac{p_{11} + p_{12}}{2}\right)$$

in seasons 1 and 2, respectively. The symmetry is flipped because the 2's have a different growing season. By applying Theorem 1, we finally obtain

COROLLARY 3. *Suppose that $\beta_{11} > \beta_{12}$ and $\beta_{22} > \beta_{21}$, and that*

$$\frac{1}{2}\left[\beta_{11}\left(1 - \frac{p_{21} + p_{22}}{2}\right) + \beta_{12}(1 - p_{21})\right] > \delta_1,$$

$$\frac{1}{2}\left[\beta_{21}(1 - p_{12}) + \beta_{22}\left(1 - \frac{p_{11} + p_{12}}{2}\right)\right] > \delta_2,$$

*where $p_{ij}$ is the density of species $i$ in the absence of the other species at the beginning of season $j$ in the mean-field equilibrium. Then, when $L$ is large, coexistence occurs.*

In the previous result, each species does better than the other in one season, but this is not necessary for coexistence. For a concrete example, suppose that

$$D = 1, \qquad \beta_{21} = 5.2, \qquad \beta_{22} = 1.0 \quad \text{and} \quad \delta_{21} = \delta_{22} = 2.0.$$

Solving the equations numerically shows that the integral of the density $\bar{u}_2(t)$ over the first season is approximately 0.366066. Thus, if we pick

$$\beta_{11} = 10000, \qquad \beta_{12} = 0, \qquad \delta_{11} = 6000 \quad \text{and} \quad \delta_{12} = 100,$$

then we have $\beta_{11}/\delta_{11} = 5/3 < 13/5 = \beta_{21}/\delta_{21}$ but

$$\frac{1}{2D}\int_0^{2D} \beta^1(t)(1 - \bar{u}_2(t))\, dt = \frac{\beta_{11}}{2}\int_0^1 (1 - \bar{u}_2(t))\, dt > 3169 > 3050$$

$$= \frac{\delta_{11} + \delta_{12}}{2},$$

so the 1's can invade the 2's. In other respects,

$$\bar{u}_1(t) \le 1 - \delta_{11}/\beta_{11} = 0.4 \qquad \text{for } t \le 1,$$

$$\bar{u}_1(t) \le 0.4\exp(-100(t-1)) \qquad \text{for } 1 \le t \le 2,$$

which implies that

$$\frac{1}{2D}\int_0^{2D} \beta^2(t)(1 - \bar{u}_1(t))\, dt = \frac{\beta_{21}}{2}\int_0^1 (1 - \bar{u}_1(t))\, dt + \frac{\beta_{22}}{2}\int_1^2 (1 - \bar{u}_1(t))\, dt$$

$$\ge \frac{3}{10}\beta_{21} + \frac{\beta_{22}}{2}\left(1 - \frac{2}{5}\int_0^1 \exp(-100t)\, dt\right)$$

$$\approx \frac{3}{10}\beta_{21} + \frac{498}{1000}\beta_{22} = 2.058 > \frac{\delta_{21} + \delta_{22}}{2},$$

so the 2's can invade the 1's as well. The fact that coexistence may occur even though species 2 is a better competitor than species 1 in both seasons is due to the combination of the presence of seasons (temporal heterogeneity) and the fact that species 1 evolves at a much faster rate than species 2. In our example, the density of 1's immediately goes to the carrying capacity in season 1 and then drop quickly to 0 in season 2. Individuals of species 2, on the contrary, evolve slowly so they never really get to full strength: Their effective average density is less than in the homogeneous case, which allows the 1's to invade the 2's, even thought the 2's are superior competitors.

**2. Mean-field model.** In this section, we prove a coexistence result for the mean-field model which is analogous to Theorem 1. A common way for showing coexistence of two species models is to prove that both species are mutually invadable, that is, species 1 can invade a community of 2's when the 1's are small in number, and vice versa. This implies that their densities are bounded away from the axes on the $u_1$–$u_2$ plane. In Proposition 2.1, below, coexistence for the mean-field model means that starting with a positive density of type 1 and type 2,

$$\liminf_{t \to \infty} u_i(t) > 0 \qquad \text{for } i = 1, 2.$$

Note that, similarly to the interacting particle system, due to temporal heterogeneities, the mean-field model (2) has no fixed point but instead converges to a limit cycle.

PROPOSITION 2.1 (Coexistence). *Suppose that*

$$\frac{1}{2D} \int_0^{2D} \beta^1(t)(1 - \bar{u}_2(t)) \, dt > \delta_1 \quad \text{and} \quad \frac{1}{2D} \int_0^{2D} \beta^2(t)(1 - \bar{u}_1(t)) \, dt > \delta_2.$$

*Then there is coexistence in the mean-field model (2).*

PROOF. By symmetry, we only need to prove survival of the 1's. The idea is to show the existence of a small $\varepsilon > 0$, such that if $u_1(t) \in (0, \varepsilon)$ for all $t \in [0, 2D]$ and

$$\int_0^{2D} \beta^1(t)(1 - \bar{u}_2(t)) - \delta_1 \, dt > 0,$$

then $u_1(2D) > u_1(0)$. Assuming that $u_1(t) \in (0, \varepsilon)$, we obtain

$$\frac{du_1}{dt} > \beta^1(t)u_1(1 - \varepsilon - u_2) - \delta_1 u_1.$$

Dividing both sides by $u_1$ and integrating from 0 to $2D$ gives

$$\log u_1(2D) - \log u_1(0) > (1 - \varepsilon)(\beta_{11} + \beta_{12})D - \int_0^{2D} \beta^1(t)u_2(t) \, dt - 2\delta_1 D.$$

Therefore, if $u_2$ satisfies the condition

$$(4) \qquad (\beta_{11} + \beta_{12})D - \int_0^{2D} \beta^1(t)u_2(t)\,dt > 2\delta_1 D$$

one can choose $\varepsilon > 0$ small so that $u_1(2D) > u_1(0)$. If $u_2(0) \leq p_{21}$ then $u_2(t) \leq \bar{u}_2(t)$ at all times, which implies that

$$(\beta_{11} + \beta_{12})D - \int_0^{2D} \beta^1(t)u_2(t)\,dt \geq \int_0^{2D} \beta^1(t)(1 - \bar{u}_2(t))\,dt > 2\delta_1 D,$$

so (4) is satisfied. This proves the result when $u_2(0) \leq p_{21}$. If on the contrary $u_2(0) > p_{21}$, then, by monotonicity of the solution with respect to the initial condition, survival of type 2 holds as well. This completes the proof. $\square$

Note that taking $\bar{u}_1(t) = 0$ in Proposition 2.1 gives the following sufficient condition for the survival of species 2 in the absence of species 1.

COROLLARY 4 (Survival of a single species). *In the absence of 1's,*

$$\liminf_{t \to \infty} u_2(t) > 0 \qquad \text{whenever } \frac{\beta_{21} + \beta_{22}}{2} > \delta_2.$$

**3. Preliminaries on periodic contact processes.** From now on, we call a *periodic contact process* a contact process on $\mathbb{Z}^2/L$ in which particles give birth at rate $\alpha(t)$ and die at rate $\delta$, where $\alpha(t)$ is positive, piecewise continuous and periodic. A particle produced at site $x \in \mathbb{Z}^2/L$ is sent to a site chosen uniformly from $\mathcal{N}(x)$. The birth occurs if the target site is empty, and is suppressed otherwise. This process is denoted by $A_t(\alpha, \delta)$ and thought of as a set valued process, that is, $A_t(\alpha, \delta)$ is the subset of $\mathbb{Z}^2/L$ corresponding to the set of sites which are occupied at time $t$. The dual process of the periodic contact process starting at site $x$ at time $t$ is denoted by $\hat{A}_s^{x,t}$ and represents the descendants of $(x, t)$ at time $t - s$. Note that contrary to the basic contact process, the periodic contact process is not self-dual since particles in the dual process $\hat{A}_s^{x,t}$ give birth at rate $\alpha(t - s)$, which a priori is different from $\alpha(s)$.

Similarly, we call a *periodic branching random walk* on $\mathbb{R}^2$ a branching random walk in which particles give birth at rate $\alpha(t)$ and die at rate $\delta$. A particle produced at $x \in \mathbb{R}^2$ is sent to a site chosen uniformly at random from $x + [-1, 1]^2$. This process is denoted by $Z_t(\alpha, \delta)$. For both the periodic contact process and the periodic branching random walk, we use a superscript to indicate the initial configuration. We first prove some convergence properties of $A_t(\alpha, \delta)$.

LEMMA 3.1 (Convergence to branching random walk). *Let $A = \{0\}$ and $T > 0$. Then, as the range $L \to \infty$, the periodic contact process $A_t^A(\alpha, \delta)$* *converges weakly to $Z_t^A(\alpha, \delta)$ in the Skorohod space of càdlàg functions from $[0, T]$ to the subsets of $\mathbb{R}^2$.*

PROOF.    This is a classical result when the birth rate $\alpha(t)$ is a constant [see, e.g., Durrett ([1991](#)) for a proof]. Let $b = \max \alpha(t)$. Then, for any $t \leq T$, the number of births by time $t$ in the contact process is dominated by the number of particles in $Z_t^A(b, 0)$ from which it follows that

$$E|A_t^A(\alpha, \delta)| \leq E|Z_t^A(b, 0)| = e^{bt} \leq e^{bT}.$$

Then, by Markov's inequality, we have

$$P(|A_t^A(\alpha, \delta)| > L^{1/3}) \leq \frac{e^{bT}}{L^{1/3}} \to 0 \qquad \text{as } L \to \infty.$$

Now, conditional on $E_1 = \{|A_t^A(\alpha, \delta)| \leq L^{1/3} \text{ for all } t \leq T\}$, at each birth, the probability that the offspring is sent to a site already occupied is less than

$$L^{1/3} \cdot \frac{L^{1/3}}{(2L+1)^2} \to 0 \qquad \text{as } L \to \infty$$

since the size of the interaction neighborhood is bounded by $(2L+1)^2$. We denote by $E_2$ the event that, up to time $T$, all the offspring of the periodic contact process are sent to empty sites. Then, on the event $E_1 \cap E_2$, one can create a coupling in which all the particles in $A_t^A(\alpha, \delta)$ are within distance $L^{1/3}/L$ in the uniform norm of their counterparts in the branching random walk. For both processes, we call particle 0 the particle at site 0 at time 0, and particle $i$, $i \geq 1$, the particle produced at the $i$th birth event. For $x \in \mathbb{R}^2$, let $\pi_L(x)$ be the (unique) vector such that

$$\|x - \pi_L(x)\|_\infty < 1/2L \quad \text{and} \quad \pi_L(x) \in \mathbb{Z}^2/L.$$

The processes are coupled as follows. If particle $i$ in the periodic branching random walk dies at time $t_i$, then particle $i$ in the periodic contact process dies at time $t_i$ as well. Now, assume that particle $i$ in the periodic branching random walk is located at $x$ and that particle $i$ in the periodic contact process is located at $x_L$. If particle $i$ in the periodic branching random walk produces an offspring which is sent to $y$ at time $t_i$, then particle $i$ in the periodic contact process produces an offspring sent to $x_L + \pi_L(y - x)$ at time $t_i$. Since there are at most $L^{1/3}$ birth events by time $T$, adding the errors we find that, on the event $E_1 \cap E_2$, the distance between particles in the contact process and their counterparts in the branching random walk is bounded by

$$L^{1/3} \sup_x \pi_L(x) = L^{1/3}/2L \to 0 \qquad \text{as } L \to \infty.$$

Since $P(E_1 \cap E_2) \to 1$ as $L \to \infty$, we conclude that the periodic contact process $A_t^A(\alpha, \delta)$ converges weakly to the periodic branching random walk $Z_t^A(\alpha, \delta)$ as $L \to \infty$. $\square$

We now prove that the number of particles in the periodic contact process $A_t(\alpha, \delta)$ is almost deterministic when $L$ is large: If we focus our attention to a space–time box of finite size, one can show that the number of particles does not deviate much from its expected value.

LEMMA 3.2 (Convergence to expected value). *Assume that $\alpha(t)$ is piecewise constant and periodic. For $A \subset \mathbb{R}^2$, let $A_t^A$ denote the periodic contact process starting from $A \cap (\mathbb{Z}^2/L)$. For any subset $A$, any site $x$ on the rescaled lattice $\mathbb{Z}^2/L$ and any time $t \geq 0$, let*

$$N(x,t) = |A_t^A \cap \{x + [0,1]^2\}|.$$

*Fix $T > 0$ and $S > 0$. Then for all $\varepsilon > 0$,*

$$P(|N(x,t) - EN(x,t)| > 5L^2\varepsilon \text{ for some}$$
$$x \in [-\sqrt{T}, \sqrt{T}]^2 \text{ and some } t \leq S) \to 0 \qquad \text{as } L \to \infty.$$

PROOF. We follow the proof of Lemma 3.5 in Durrett and Lanchier (2008), in which the result is proved when the birth rate $\alpha(t)$ is a constant. Let

$$N^\varepsilon(x,t) = |A_t^A \cap \{x + [0,\varepsilon]^2\}| \quad \text{and} \quad y,z \in x + [0,\varepsilon]^2 \cap (\mathbb{Z}^2/L).$$

First, we observe that the number of descendants in the dual processes $\hat{A}_t^{y,t}$ and $\hat{A}_t^{z,t}$ depends on $t$ and the birth and death rates, but not on $L$. In particular, since the interaction neighborhood contains $(2L+1)^2 - 1$ sites, the probability that the dual processes hit, that is, that a descendant is sent to a site already occupied by another descendant, is $\leq C_1/L^2$. Now, noting that

$$A_t^A(y) = \begin{cases} 0, & \text{if } \hat{A}_t^{y,t} \cap A = \varnothing, \\ 1, & \text{if } \hat{A}_t^{y,t} \cap A \neq \varnothing \end{cases} \quad \text{and} \quad A_t^A(z) = \begin{cases} 0, & \text{if } \hat{A}_t^{z,t} \cap A = \varnothing, \\ 1, & \text{if } \hat{A}_t^{z,t} \cap A \neq \varnothing, \end{cases}$$

a standard construction shows that the covariance of $A_t^A(y)$ and $A_t^A(z)$ can be bounded by the probability that the dual processes hit [see page 21 in Griffeath (1979)]:

$$(5) \qquad \qquad \text{cov}(A_t^A(y), A_t^A(z)) \leq C_1/L^2.$$

Since $N^\varepsilon(x,t) = \sum_{y \in x + [0,\varepsilon]^2} A_t^A(y)$, inequality (5) implies that

$$\text{Var } N^\varepsilon(x,t) \leq \varepsilon^2 L^2 + (\varepsilon^2 L^2)^2 C_1/L^2 \leq C_2 L^2 \varepsilon^2.$$

Then, by Chebyshev's inequality, we obtain

(6)     $P(|N^\varepsilon(x,t) - EN^\varepsilon(x,t)| > L^2\varepsilon^3) \leq C_2/\varepsilon^4 L^2 \to 0$     as $L \to \infty$.

This holds for all $x \in [-\sqrt{T}, \sqrt{T} + 1]^2 \cap (\mathbb{Z}^2/L)$ and all $t \leq S$. To deduce the result, we first use a discretization of space and time. The idea is to rely on the convergence in (6) to control the densities $N^\varepsilon(x,t)$ at a finite number (that does not depend on $L$) of space–time points. More precisely, we let $\tau > 0$ be a positive constant and look at the densities $N^\varepsilon(x,t)$ for space–time points $(x,t)$ belonging to a subset of $\varepsilon\mathbb{Z}^2 \times \tau\mathbb{Z}_+$. As $L \to \infty$,

$P(|N^\varepsilon(x,n\tau) - EN^\varepsilon(x,n\tau)| > L^2\varepsilon^3$ for some $0 \leq n \leq m$ and

$\qquad\qquad$ some $x \in [-\sqrt{T}, \sqrt{T} + 1]^2 \cap \varepsilon\mathbb{Z}^2)$

$\qquad \leq C_3 T(m+1)/\varepsilon^6 L^2 \to 0,$

where $m = \min\{n \geq 1 : n\tau \geq S\}$. Using large deviation estimates for the Binomial distribution and taking $\tau > 0$ small, one can show that the number of sites in $x + [0,\varepsilon]^2$ that flip between $n\tau$ and $(n+1)\tau$ is smaller than $2L^2\varepsilon^3$ with probability at least $1 - C_3 \exp(-\alpha L^2)$ so that

$P(|N^\varepsilon(x,t) - EN^\varepsilon(x,t)| > 3L^2\varepsilon^3$ for some $t \leq S$

$\qquad\qquad$ and some $x \in [-\sqrt{T}, \sqrt{T} + 1]^2 \cap \varepsilon\mathbb{Z}^2) \to 0$     as $L \to \infty$.

To conclude, we observe that if

$|N^\varepsilon(x,t) - EN^\varepsilon(x,t)| \leq 3L^2\varepsilon^3$     for all $x \in [-\sqrt{T}, \sqrt{T} + 1]^2 \cap \varepsilon\mathbb{Z}^2$,

then we have that

$|N(x,t) - EN(x,t)| \leq 3L^2\varepsilon^3/\varepsilon^2 + 2L^2\varepsilon = 5L^2\varepsilon$

for all $x \in [-\sqrt{T}, \sqrt{T}]^2 \cap (\mathbb{Z}^2/L)$ since there are at most $\varepsilon^{-2}$ squares with length side $\varepsilon$ included in the unit square $x + [0,1]^2$ and at most $2L^2\varepsilon$ remaining sites. $\quad\square$

LEMMA 3.3 (Convergence to ODE). *Let $A_t = A_t(\beta^1, \delta)$ be the periodic contact process starting from the "all occupied" configuration. For any site $x$ on the rescaled lattice $\mathbb{Z}^2/L$, we set*

$$u_L(x,t) = P(x \in A_t(\beta^1, \delta)).$$

*Let $T > 0$. Then, as the range of the interactions $L \to \infty$, the function $t \mapsto u_L(x,t)$ converges uniformly on $[0,T]$ to a solution $u(t)$ of the ODE*

$$\frac{du}{dt} = \beta^1(t)u(1-u) - \delta u \qquad \text{with } u(0) = 1.$$

PROOF. Since the evolution rules are translation invariant, the probability of $\{x \in A_t\}$ does not depend on $x$. Therefore, we can set $u_L(x,t) = u_L(t)$. Let $T > 0$ and $S(t)$ be the semigroup that generates the process $A_t$. For any cylindric function $f$, that is any function that depend on finitely many coordinates, the Hille–Yosida theorem states that

$$(7) \qquad \frac{d}{dt} S(t)f = S(t)\Omega f,$$

where $\Omega$ is the Markov generator of $S(t)$. Now, letting $f(A_t) = \mathbb{1}_{\{x \in A_t\}}$, we have

$$(8) \qquad S(t)f = Ef(A_t) = P(x \in A_t) = u_L(t).$$

Combining (7) and (8) implies that

$$(9) \qquad \frac{d}{dt} P(x \in A_t) = \sum_{y \in \mathcal{N}(x)} \frac{\beta^1(t)}{|\mathcal{N}(x)|} P(x \notin A_t \text{ and } y \in A_t) - \delta P(x \in A_t).$$

We claim that $\{x \in A_t\}$ and $\{y \in A_t\}$ are asymptotically independent for any $t \le T$, that is,

$$(10) \qquad P(x \notin A_t \text{ and } y \in A_t) \to P(x \notin A_t) \cdot P(y \in A_t) \qquad \text{as } L \to \infty.$$

To see this, we first observe that (10) is equivalent to

$$(11) \qquad P(\hat{A}_t^{x,t} = \varnothing \text{ and } \hat{A}_t^{y,t} \ne \varnothing) \to P(\hat{A}_t^{x,t} = \varnothing) \cdot P(\hat{A}_t^{y,t} \ne \varnothing) \qquad \text{as } L \to \infty.$$

Since $t \le T$, the number of particles in $\hat{A}_s^{x,t}$ and $\hat{A}_s^{y,t}$, $s \le t$, depends on $T$ and the birth and death rates but not on $L$ so, with probability arbitrarily close to 1 when $L$ is sufficiently large, the number of particles is bounded by $L^{1/3}$ (see the proof of Lemma 3.1). Moreover, on such an event, the probability that the dual processes hit is less than $2L^{1/3}/4L^2 \to 0$. In particular, the dual processes are independent in the limit as $L \to \infty$ which establishes (11) and (10). Finally, combining (9) and (10) gives the desired result. $\square$

## 4. Proof of Theorem 1.

4.1. *Block construction.* The first step to prove Theorem 1 is to construct the periodic competition model $\xi_t$ from a Harris' graphical representation [Harris (1972)]. For $x, y \in \mathbb{Z}^2/L$, we let

$$q(x,y) = \begin{cases} |\mathcal{N}(x)|^{-1}, & \text{if } 0 < \|x - y\|_\infty \le 1, \\ 0, & \text{otherwise.} \end{cases}$$

For any ordered pair of sites $(x,y)$ and any $i,j \in \{1,2\}$, we introduce the following collections of independent Poisson processes:

| arrival times | rates | symbols |
|---|---|---|
| $T_n^{ij}(x,y)$ | $\beta_{ij} q(x,y)$ | $x \xrightarrow{i} y$ |
| $U_n^i(x)$ | $\delta_i$ | $\times$ at $x$ |

that is, we draw an $i$-arrow from site $x$ to site $y$ at time $T_n^{ij}(x, y)$ to indicate that if $x$ is occupied by species $i$, site $y$ is empty, and

$$(12) \qquad \lceil T_n^{ij}(x, y)/D \rceil := \min\{a \in \mathbb{Z} : a \geq T_n^{ij}(x, y)/D\} \equiv j \pmod{2}$$

then site $y$ becomes occupied by an individual of species $i$. Note that the condition in (12) means that the arrival time $T_n^{ij}(x, y)$ occurs in season $j$. In addition, we put a $\times$ at site $x$ at time $U_n^i(x)$ to indicate that an individual of species $i$ at $x$ dies. A result of Harris (1972) implies that the previous construction induces a well-defined Markov process. Moreover, the flip rates indicate that this process is the periodic competition model.

Our proof of Theorem 1 relies on the use of a block construction, which is now a standard technique. The idea is to show that, when viewed on suitable length and time scales, our process dominates the set of wet sites of a $M$-dependent oriented percolation process on

$$\mathcal{L} = \{(m, n) \in \mathbb{Z}^2 : m + n \text{ is even and } n \geq 0\}$$

in which sites are open with probability $p$ close to 1. More precisely, each site $(m, n) \in \mathcal{L}$ is associated with a random variable $\omega(m, n) \in \{0, 1\}$ which indicates whether the site is open (state 1) or closed (state 0). The $M$-dependency means that

$$P(\omega(m_i, n_i) = 1 \text{ for } 1 \leq i \leq k) = p^k,$$

whenever $\|(m_i, n_i) - (m_j, n_j)\| > M$ for $i \neq j$. A site $(m, n)$ is said to be wet if there exists a sequence of integers $m_0, m_1, \ldots, m_n = m$ such that:

1. For $i = 0, 1, \ldots, n - 1$, we have $\|m_{i+1} - m_i\| = 1$.
2. For $i = 0, 1, \ldots, n$, site $(m_i, i)$ is open, that is $\omega(m_i, i) = 1$.

See Durrett (1995) for more details. To compare the periodic competition model viewed on suitable length and time scales with oriented percolation, we let $\gamma$ be a large parameter to be fixed later, and let $N_i(x, t)$ denote the number of type $i$ particles in $x + [0, 1]^2$ at time $t$. Assuming that

$$(13) \qquad N_i(0, 0) > L^2 \exp(-\gamma T) \qquad \text{for } i = 1, 2,$$

we will prove that with probability close to 1 when the parameters $T$ and $L$ are large,

$$(14) \qquad N_i((\sqrt{T}, 0), T) > L^2 \exp(-\gamma T) \qquad \text{for } i = 1, 2.$$

In particular, if we call $(m, n) \in \mathcal{L}$ an occupied site whenever

$$N_i(m(\sqrt{T}, 0), nT) > L^2 \exp(-\gamma T) \qquad \text{for } i = 1, 2,$$

then the set of occupied sites dominates the set of wet sites of an oriented percolation process in which sites are open with probability $p$ close to 1.

The construction of the percolation process relies on a coupling argument in which the random variables $\omega(m, n)$ are defined by induction over the level $n$. This technique is now standard and we refer the reader to the Appendix of Durrett (1995) for a detailed construction of the process in the general case. Since percolation occurs whenever the density $p$ is sufficiently close to 1 (i.e., starting with one open site at level 0, there is a positive probability that the cluster of wet sites is infinite), Theorem 1 follows from the comparison with oriented percolation. In conclusion, the aim is to prove that the conditional probability of the events in (14) given the events in (13) can be made arbitrarily close to 1.

Since the proofs of the inequalities in (14) for $i = 1$ and $i = 2$ are identical, we will only show the first one ($i = 1$). There are three steps in showing the survival of species 1 (see Figure 2 for an illustration of the three steps described below):

(i) The density of type 2 particles in $[-\sqrt{T}, \sqrt{T}]^2$ between times $\sqrt{T}$ and $T$ is (almost) bounded by the mean-field equilibrium $\bar{u}_2(t)$.

(ii) If the unit square $[0, 1]^2$ contains at least $L^2 \exp(-\gamma T)$ type 1 particles at time 0, then it will contain at least $L^2 \exp(-\gamma T + cT)$ type 1 particles at time $T - 2\sqrt{T}$ for some $c > 0$.

(iii) If the unit square $[0, 1]^2$ contains at least $L^2 \exp(-\gamma T + cT)$ type 1 particles at $T - 2\sqrt{T}$, then at least $L^2 \exp(-\gamma T)$ of them will "invade" the unit square $(\sqrt{T}, 0) + [0, 1]^2$ at time $T$.

4.2. *Proofs of (i)–(iii).* We start with condition (i) which is proved in Lemma 4.1.

LEMMA 4.1. *For any $\varepsilon > 0$ and for $T$ sufficiently large,*

$$P(N_2(x, t) > L^2(\bar{u}_2(t) + 7\varepsilon) \text{ for some}$$

$$x \in [-\sqrt{T}, \sqrt{T}]^2 \text{ and some } t \in [\sqrt{T}, T]) \to 0 \qquad \text{as } L \to \infty.$$

PROOF. Let $A_t$ be the periodic contact process with birth rate $\beta^2(t)$, death rate $\delta_2$, and initial configuration $A_0 \equiv 2$, and let $\zeta_t = \{x : \xi_t(x) = 2\}$ denote the set of sites occupied by a type 2 particle at time $t$ in the periodic competition model. Then $A_t$ and $\xi_t$ can be coupled in such a way that $\zeta_t \subset A_t$ for all $t \geq 0$. In particular, it suffices to prove the result for

$$N(x, t) = |A_t \cap \{x + [0, 1]^2\}|$$

instead of $N_2(x, t)$. Let $u(t)$ be the solution of the ODE

$$\frac{du}{dt} = \beta^2(t) u(1 - u) - \delta_2 u \qquad \text{with } u(0) = 1.$$

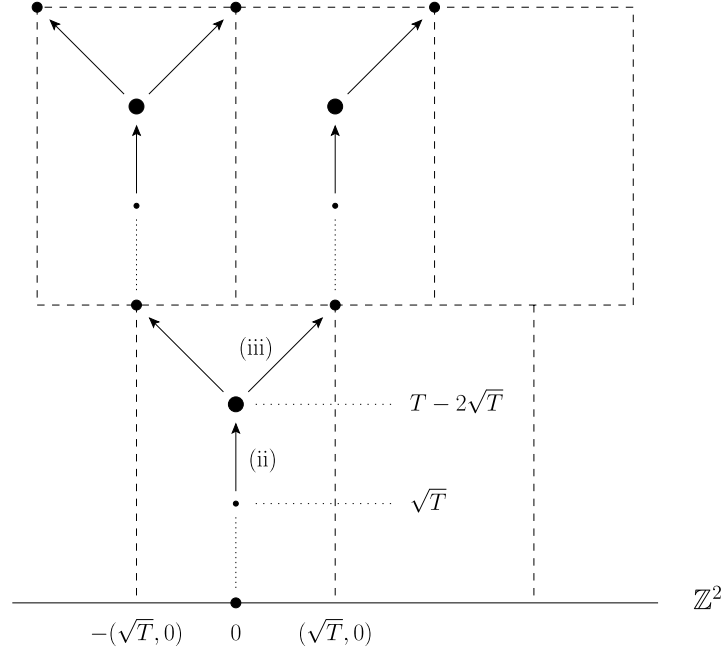

Fig. 2.   *Boxes in the block construction for Theorem 1.*

Then there exists $T$ sufficiently large such that

(15)             $|u(t) - \bar{u}_2(t)| < \varepsilon$        for all $t \geq \sqrt{T}$.

Moreover, by Lemma 3.3, $P(x \in A_t) \to u(t)$ uniformly on $[0, T]$ as $L \to \infty$ which, together with inequality (15), implies that

(16)
$$EN(x, t) < L^2(\bar{u}_2(t) + 2\varepsilon)$$

for all $(x, t) \in [-\sqrt{T}, \sqrt{T}]^2 \times [\sqrt{T}, T]$ as $L \to \infty$.

In other respects, by Lemma 3.2,

(17)
$$P(|N(x, t) - EN(x, t)| > 5L^2\varepsilon \text{ for some}$$

$x \in [-\sqrt{T}, \sqrt{T}]^2$ and some $t \in [\sqrt{T}, T]) \to 0$     as $L \to \infty$.

The result follows by combining (16) and (17).   □

Let $\bar{Z}_t = \bar{Z}_t(\alpha, \delta)$ be the restriction of the branching random walk with birth rate $\alpha(t)$ and death rate $\delta$ to the square $I_T = [-\sqrt{T}, \sqrt{T}]^2$, that is, particles landed outside $I_T$ are killed. In preparation for proving (ii), we need to show that if the averaged birth rate $\alpha(t)$ along one period is greater than the death rate $\delta$ then the number of particles in $\bar{Z}_t$ grows exponentially fast.

LEMMA 4.2. *Assume that $\bar{Z}_0 = \bar{Z}_0^x = \{x\}$. For $s > r > 0$ fixed and $T$ sufficiently large*

$$E|\bar{Z}_t^x \cap \{y + [0,1]^2\}| \geq \exp\left(\frac{1}{2}\int_0^t (\alpha(\theta) - \delta)\, d\theta\right)$$

*for all $x, y \in [-\sqrt{T}/2, \sqrt{T}/2]^2$ and $t \in [rT, sT]$.*

PROOF. Let $S_t^x = S_t^x(\alpha)$ be the random walk that jumps at rate $\alpha(t)$ and has kernel the uniform distribution on the square $[-1,1]^2$, that is, the new position of the process is chosen uniformly at random from $S_t^x + [-1,1]^2$, and denote by $\bar{S}_t^x$ its restriction to the square $I_T$. Let

$$m(t, x, A) = E|\bar{Z}_t^x(\alpha, \delta) \cap A|$$

denote the mean number of particles in the set $A$ at time $t$. Then we claim that

$$(18) \qquad m(t, x, A) = \exp\left(\int_0^t (\alpha(\theta) - \delta)\, d\theta\right) P(\bar{S}_t^x(\alpha) \in A).$$

This holds because both sides of (18) satisfy the differential equation

$$\frac{\partial m(t, x, A)}{\partial t} = -\delta m(t, x, A) + \int \alpha(t) m(t, x, dy) \mathbb{1}_{\{y \in [-\sqrt{T}, \sqrt{T}]^2\}} \nu(A - y),$$

where $\nu$ is the uniform probability measure on $[-1,1]^2$. Since $0 < \min \alpha(t) \leq \max \alpha(t) < \infty$, one can rescale time so that $\bar{S}_t^x$ has constant jump rate $\min \alpha(t)$. Let

$$t = uF(T), \qquad x = v\sqrt{F(T)}, \qquad y = w\sqrt{F(T)}, \qquad v, w \in [-\tfrac{1}{2}, \tfrac{1}{2}]^2.$$

Here $F(T) = \mathcal{O}(T)$ is piecewise continuous with $1 \leq F(T)/T \leq \max \alpha(t)$. An application of the invariance principle shows that the process $\bar{S}_t^x$ converges in probability to a Brownian motion when time and space are properly scaled. More precisely, sending $T \to \infty$ we have

$$F(T) \cdot P(\bar{S}_t^x \in y + [0,1]^2) \to \bar{p}_u(v, w),$$

where $\bar{p}_u(v, w)$ is the time-$u$ transition probability of a 2-dimensional Brownian motion from $v$ to $w$ that is killed outside $[-1,1]^2$. Since $|v - w| \leq 1$, the probability $\bar{p}_u(v, w)$ is bounded from below which implies that, for $T$ sufficiently large,

$$P(\bar{S}_t^x \in y + [0,1]^2) = \mathcal{O}(T^{-1}).$$

It follows that

$$m(t, x, y + [0,1]^2) = \exp\left(\int_0^t (\alpha(\theta) - \delta)\, d\theta\right) \cdot \mathcal{O}(T^{-1}) \geq \exp\left(\frac{1}{2}\int_0^t (\alpha(\theta) - \delta)\, d\theta\right)$$

for all $t \in [rT, sT]$ with $T$ sufficiently large.  $\square$

In the next lemma, we establish (ii). The idea is that, since the density of type 2 particles is bounded by $\bar{u}_2 + 7\varepsilon$, the 1's grow exponentially fast from time $\sqrt{T}$ to time $T - 2\sqrt{T}$ so that they can regain all possible losses occurring by time $\sqrt{T}$.

LEMMA 4.3.  *Assume that $N_1(0,0) > L^2 \exp(-\gamma T)$. Then, for $\gamma$ and $T$ sufficiently large,*

$$P(N_1(0, T - 2\sqrt{T}) \le L^2 \exp(-\gamma T + cT)) \to 0 \qquad as \ L \to \infty$$

*for a suitable constant $c > 0$.*

PROOF.  First of all, we observe that the probability that a type 1 particle remains alive $\sqrt{T}$ units of time is given by $p = \exp(-\delta_1\sqrt{T})$. Therefore, an application of the large deviation inequality for the Binomial distribution implies that as $L \to \infty$,

$$P(N_1(0, \sqrt{T}) \le (L^2/2) \exp(-\gamma T - \delta_1\sqrt{T})) \le \exp(-e^{-\gamma T - \delta_1\sqrt{T}} L^2/8) \to 0.$$

In particular, since the inclusion relation $\subset$ on the set of sites occupied by type 1 particles is preserved by the dynamics, it suffices to prove the result when

(19)
$$N_1(0, \sqrt{T}) = (L^2/2) \exp(-\gamma T - \delta_1\sqrt{T}) \quad \text{and}$$
$$\xi_{\sqrt{T}}(x) \ne 1 \qquad \text{for all } x \notin [0,1]^2.$$

The assumptions of Theorem 1 allow to fix $\varepsilon > 0$ such that

(20)
$$\int_0^{2D} \beta^1(t)(1 - \bar{u}_2(t) - 8\varepsilon) - \delta_1 \, dt > 0.$$

Applying (18) to the restricted branching random walk $\bar{Z}_t$ with birth rate $\alpha(t) = \beta^1(t)$ and death rate $\delta = \delta_1$, which dominates the number of type 1 particles in the square $I_T$, and taking (19) into account, we obtain that, for all sites $x \in I_T$ and all times $t \in [\sqrt{T}, T]$,

$$EN_1(x,t) \le N_1(0, \sqrt{T}) \exp\left(\int_{\sqrt{T}}^t (\beta^1(\theta) - \delta_1) \, d\theta\right)$$
$$\le L^2 \exp(-\gamma T + \max(\beta_{11}, \beta_{12})T).$$

In particular, there exists a sufficiently large $\gamma$, fixed from now on, such that

(21)  $P(N_1(x,t) \le L^2\varepsilon \text{ for all } x \in I_T \text{ and all } t \in [\sqrt{T}, T - 2\sqrt{T}]) \to 1$

$$as \ L \to \infty.$$

By Lemma 4.1, the density of type 2 particles inside $I_T = [-\sqrt{T}, \sqrt{T}]^2$ is bounded by $\bar{u}_2 + 7\varepsilon$ between time $\sqrt{T}$ and time $T - 2\sqrt{T}$. This, together with (21), implies that, with probability arbitrarily close to 1 when $L$ is large, the number of empty sites in each square $x + [0,1]^2$ is larger than $L^2(1 - \bar{u}_2(t) - 8\varepsilon)$. So, in a given interval, as $L \to \infty$, the number of 1's dominates the number of particles in a branching random walk with birth rate

$$\alpha(t) = \beta^1(t)(1 - \bar{u}_2(t) - 8\varepsilon) \tag{22}$$

and death rate $\delta_1$, and with $(L^2/2)\exp(-\gamma T - \delta_1\sqrt{T})$ particles in $[0,1]^2$ at time $\sqrt{T}$. To see this, note first that (22) is a lower bound of the rate of successful birth of type 1, that is, the rate at which type 1 particles are actually created, in the periodic contact process when (21) holds. Moreover, the arguments of Lemma 3.1 imply that both processes can be coupled in such a way that type 1 particles in the periodic contact process are arbitrarily close to their counterparts in the branching random walk. In particular, it suffices to prove the result for the branching random walk described above restricted to the square $I_T$. We denote this process by $\bar{Z}_t(\alpha, \delta_1)$ and let

$$\bar{N}(x,t) = |\bar{Z}_t(\alpha, \delta_1) \cap \{x + [0,1]^2\}|.$$

By Lemma 4.2 and (20), for $T$ sufficiently large,

$$E\bar{N}(0, T - 2\sqrt{T}) > (L^2/2)\exp(-\gamma T - \delta_1\sqrt{T})$$
$$\times \exp\left(\frac{1}{2}\int_{\sqrt{T}}^{T-2\sqrt{T}} \beta^1(t)(1 - \bar{u}_2(t) - 8\varepsilon) - \delta_1\, dt\right)$$
$$> L^2\exp(-\gamma T + 2cT)$$

for a suitable constant $c > 0$. Finally, by Lemma 3.2 (which applies to the periodic branching random walk as well),

$$\bar{N}(0, T - 2\sqrt{T}) > L^2\exp(-\gamma T + 2cT) - 5L^2\varepsilon > L^2\exp(-\gamma T + cT)$$

with probability arbitrarily close to 1 for $L$ and $T$ sufficiently large. □

To complete the proof of Theorem 1, it remains to "move" our particles from the square $[0,1]^2$ to the square $(\sqrt{T}, 0) + [0,1]^2$ and show that $L^2\exp(-\gamma T)$ of them will live until time $T$. Our last lemma shows (iii) in the list.

LEMMA 4.4. *Assume that* $N_1(0, T - 2\sqrt{T}) > L^2\exp(-\gamma T + cT)$ *and that*

$N_2(x,t) \leq L^2(\bar{u}_2(t) + 7\varepsilon)$ *for all* $x \in [-\sqrt{T}, \sqrt{T}]^2$ *and* $t \in [T - 2\sqrt{T}, T]$.

*Then, for $T$ large, $P(N_1((\sqrt{T}, 0), T) \leq L^2\exp(-\gamma T)) \to 0$ as $L \to \infty$.*

PROOF.    Fix $a > 0$. For any integer $m \in \mathbb{Z}$ and any time $t \geq 0$, we let

$$B_m = (m/2, 0) + [0, 1/2]^2 \quad \text{and} \quad H_{m,t} = |\{x \in B_m : \xi_t(x) = 1\}|.$$

Assume first that $H_{m+1,t} \geq 2L^2 \exp(-aT)$ for some $t \in (s, s+1)$. The probability that a type 1 particle remains alive at least one unit of time is $e^{-\delta_1}$. Since these events are independent, the random variable $H_{m+1,s+1}$ dominates a Binomial random variable with parameters

$$n \geq 2L^2 \exp(-aT) \quad \text{and} \quad p = e^{-\delta_1}.$$

An application of the large deviation result for the Binomial

$$P(X \leq n(p - z)) \leq \exp(-nz^2/2p)$$

with $z = p/2$ then gives

$$
\begin{aligned}
P(H_{m+1,s+1} &\leq L^2 \exp(-aT - \delta_1)|H_{m+1,t} \geq 2L^2 \exp(-aT) \\
&\qquad\qquad\qquad\qquad \text{for some } t \in (s, s+1)) \\
&\leq \exp(-np/8) \leq \exp(-e^{-aT-\delta_1}L^2/4).
\end{aligned}
$$
(23)

Now, assume on the contrary that $H_{m+1,t} < 2L^2 \exp(-aT)$ for all $t \in (s, s+1)$. Then the probability that a type 1 particle in $B_m$ at time $s$ gives birth by time $s+1$ to an offspring which is sent to an empty site in $B_{m+1}$ and that both the parent and the offspring live until time $s+1$ is bounded from below by

$$p_0 = \tfrac{1}{16} e^{-2\delta_1}(1 - e^{-b})(1 - \max(p_{21}, p_{22}) - 7\varepsilon - 2\exp(-aT)),$$

where $b = \min(\beta_{11}, \beta_{12})$. The terms on the right-hand side are lower bounds for the probabilities that (i) the offspring is sent to $B_{m+1}$, (ii) particles live until time $s+1$, (iii) the birth occurs by time $s+1$, and (iv) the landing site is empty (since $\bar{u}_2(t) \leq \max(p_{21}, p_{22})$). This implies that

$$
\begin{aligned}
P(H_{m+1,s+1} &\leq (p_0/2)H_{m,s}|H_{m+1,t} < 2L^2 \exp(-aT) \\
&\qquad \text{for all } t \in (s, s+1) \text{ and } H_{m,s} > L^2 \exp(-\gamma T)) \to 0
\end{aligned}
$$
(24)

as $L \to \infty$. Applying (23) and (24) with $m = 2\sqrt{T} - n$ and $s = T - n$ gives

$$H_{2\sqrt{T}-n+1, T-n+1} > \tfrac{1}{2} \min(p_0, e^{-\delta_1})H_{2\sqrt{T}-n, T-n} = \frac{p_0}{2} H_{2\sqrt{T}-n, T-n}$$

with probability close to 1 when $L$ is large on the event that $H_{m,s} > L^2 \exp(-\gamma T)$. Assuming that one fourth of the particles of $N_1(0, T - 2\sqrt{T})$ are found in $H_{0, T-2\sqrt{T}}$ and observing that

$$\frac{1}{4}\left(\frac{p_0}{2}\right)^{2\sqrt{T}} L^2 \exp(-\gamma T + cT) > L^2 \exp(-\gamma T)$$

for $T$ large enough, a simple induction allows us to conclude that

$$H_{2\sqrt{T},T} > \frac{p_0}{2} H_{2\sqrt{T}-1,T-1} > \left(\frac{p_0}{2}\right)^{2\sqrt{T}} H_{0,T-2\sqrt{T}}$$

$$> \frac{1}{4}\left(\frac{p_0}{2}\right)^{2\sqrt{T}} L^2 \exp(-\gamma T + cT)$$

with probability close to 1 when $L$ is large. Recalling the definition of $H_{m,t}$, we obtain that, for $T$ sufficiently large, and with probability arbitrarily close to 1 when $L$ is large,

$$N_1((\sqrt{T},0),T) > L^2 \exp(-\gamma T).$$

This completes the proof. $\square$

**Acknowledgment.** The authors would like to thank three anonymous referees for many comments that helped to improve the article.

B. CHAN                                      N. LANCHIER
R. DURRETT                                   DEPARTMENT OF MATHEMATICS
DEPARTMENT OF MATHEMATICS                       AND STATISTICS
CORNELL UNIVERSITY                           ARIZONA STATE UNIVERSITY
ITHACA, NEW YORK 14853                       TEMPE, ARIZONA 85287
USA                                          USA
                                             E-MAIL: lanchier@math.asu.edu